\documentclass[11pt,openbib]{article}

\usepackage{latexsym}
\usepackage{amsmath}
\usepackage{amssymb}
\usepackage{amsfonts}
\usepackage[breaklinks,
  colorlinks=true,
  citecolor=cyan,
  linkcolor=magenta,
  urlcolor=blue,
  hyperfootnotes=false]{hyperref}

\numberwithin{equation}{section}
\newtheorem{thrm}{Theorem}[section]

\newtheorem{defn}{Definition}[section]
\newtheorem{lemm}{Lemma}[section]

\newenvironment{rmk}{{\bf Remark:}}{\hfill{$\Box$}}

\newenvironment{prf}{{\em Proof:}}{\hfill{$\Box$}}

\newcommand{\alf}{\alpha}

\newcommand{\apx}{\approx}

\newcommand{\dt}[1]{\frac{d}{dt}{#1}}

\newcommand{\Dlt}{\Delta}

\newcommand{\dert}{\frac{d}{dt}}

\newcommand{\field}[1]{\mathbb{#1}}

\newcommand{\Gmm}{\Gamma}

\newcommand{\gs}{\geqslant}

\newcommand{\half}{\frac{1}{2}}
\newcommand{\hf}{\half}
\newcommand{\hft}{\frac{3}{2}}

\newcommand{\Lmd}{\Lambda}
\newcommand{\lan}{\langle}

\newcommand{\lf}{\left}
\newcommand{\lra}[2]{\left\langle #1,\,#2\right\rangle}
\newcommand{\ls}{\leqslant}
\newcommand{\lsp}{\preccurlyeq}

\newcommand{\nbl}{\nabla}

\newcommand{\ol}{\overline}

\newcommand{\Omg}{\Omega}
\newcommand{\ox}{\otimes}
\newcommand{\prt}{\partial}

\newcommand{\Di}{\prt_i}
\newcommand{\Pf}{\field{P}}

\newcommand{\Pn}{\prt_n}
\newcommand{\Px}{\prt_x}
\newcommand{\Py}{\prt_y}
\newcommand{\Pz}{\prt_z}

\newcommand{\qt}{\frac{1}{4}}
\newcommand{\qtt}{\frac{3}{4}}
\newcommand{\R}{\field{R}}

\newcommand{\ra}{\rightarrow}
\newcommand{\ran}{\rangle}

\newcommand{\rt}{\right}

\newcommand{\tht}{\theta}

\newcommand{\veps}{\varepsilon}

\newcommand{\vphi}{\varphi}
\newcommand{\vfi}{\varphi}

\newcommand{\x}{\times}

\begin{document}
\pagestyle{myheadings}

\title{\bf Global Uniform Boundedness of Solutions\\
 to viscous 3D Primitive Equations\\
 with Physical Boundary Conditions}

\author{ Ning Ju $^1$}
\footnotetext[1]{
	Department of Mathematics, Oklahoma State University,
	401 Mathematical Sciences, Stillwater OK 74078, USA.
	Email: {\tt ning.ju@okstate.edu},
	}
\date{}

\maketitle

\begin{abstract}
 Global uniform boundedness of solutions to 3D viscous Primitive Equations in
 a bounded cylindrical domain with {\em physical boundary condition} is proved
 in space $H^m$ for any $m\gs 2$. A bounded absorbing set for the solutions in
 $H^m$ is obtained. These results seem quite difficult to be proved using the
 methods recently developed in \cite{j:17} and \cite{j.t:15}. A completely
 different approach based on hydrostatic helmholtz decomposition is presented,
 which is also applicable to the cases with other boundary conditions. Several
 important results about hydrostatic Leray projector are obtained and utilized.
 These results are expected to be of general interest and will be helpful for
 solving some other problems for 3D viscous Primitive Equations which appeared
 hard previously for the cases with non-periodic boundary conditions (see e.g.
 \cite{j:17.g}).\\



\end{abstract}

%\tableofcontents

\markboth{3D Viscous Primitive Equations}
{\hspace{.2in} N. Ju \hspace{.5in} 3D Viscous Primitive Equations}

\indent
\baselineskip 0.58cm

\section{Introduction}
\label{s:intro}

 Let $D$ be a bounded open subset of $\R^2$ with smooth boundary $\prt D$ and
\[ \Omg =D\times(-h, 0)\subset\R^3, \]
 where $h$ is a positive constant. Consider in the cylinder $\Omg$ the system
 of viscous Primitive Equations (PEs) of Geophysical Fluid Dynamics:

\noindent
 {\em Conservation of horizontal momentum:}
\begin{equation*}
  \label{e:v}
  v_t + (v\cdot\nbl)v + w v_z +\nbl p + fv^\bot + L_1 v = 0;
\end{equation*}
 {\em Hydrostatic balance:}
\begin{equation*}
  \label{e:hs}
  p_z + \tht =0 ;
\end{equation*}
 {\em Continuity equation:}
\begin{equation*}
  \label{e:cnt}
  \nbl\cdot v + w_z = 0 ;
\end{equation*}
 {\em Heat conduction:}
\begin{equation*}
  \label{e:t}
  \tht_t + v\cdot\nbl \tht + w\tht_z + L_2 \tht=Q.
\end{equation*}

 The unknowns in the above system of 3D viscous PEs are the fluid velocity field
 $(v, w)=(v_1, v_2, w)\in\R^3 $ with $v = (v_1, v_2)$ and $v^\bot=(-v_2, v_1)$
 being horizontal, the temperature $\tht$ and the pressure $p$. The Coriolis
 rotation frequency $f = f_0(\beta + y)$ in the $\beta$-plane approximation and
 the heat source $Q$ are given. The differential operators $L_1$ and $L_2$ are
 defined respectively as:
\begin{equation*}
 \label{e:L1.L2}
   L_i := - \nu_i\Dlt - \mu_i {\Pz^2},
\end{equation*}
 with {\em positive} constants $\nu_i, \mu_i$ for $i =1, 2$. In the above
 equations, $\nbl$ and $\Dlt$ denote horizontal gradient and Laplacian:
\[ \nbl := (\Px, \Py)  \equiv (\prt_1,\prt_2), \quad
   \Dlt := \Px^2+\Py^2 \equiv \sum_{i=1}^2\Di^2. \]
 In the sequel, we also denote
\[ \nbl_3 := (\nbl, \Pz) = (\Px,\Py,\Pz). \]
 The boundary of $\Omg$ is partitioned into three parts:
 $\partial\Omg = \Gmm_t \cup \Gmm_b \cup \Gmm_l$, where
\begin{align*}
\Gmm_t &:= \{(x, y, z) \in \overline{\Omg} : z = 0\},\\
\Gmm_b &:= \{(x, y, z) \in \overline{\Omg} : z = -h\},\\
\Gmm_l &:= \{(x, y, z) \in \overline{\Omg}:
 (x,y)\in \prt D\}.
\end{align*}
 The following set of physical boundary conditions will be used:
\begin{equation*}
\begin{aligned}
  \mbox{on} \quad \Gmm_t \mbox{\,:}& \quad 
  v_z+\alf_1 v=0, \quad w=0, \quad \tht_z+\alf_2\tht=0,\\
  \mbox{on} \quad \Gmm_b \mbox{\,:}& \quad
  v=0,  \quad w=0, \quad \tht_z=0,\\
  \mbox{on} \quad \Gmm_l \mbox{\,:}& \quad 
  v =0,  \quad \Pn\tht=0,
\end{aligned}
\end{equation*}
 where $\alf_1,\alf_2$ are {\em non-negative} constants and $n$ is the normal
 vector of $\Gmm_l$.

 Let $p_s(x,y,t)$ be the pressure on the top $\Gmm_t$. Then, the above system
 of PEs can be re-written as
\begin{equation}
\label{e:v.n}
\begin{split}
 v_t + L_1v + (v\cdot\nbl)v + wv_z 
 + \nbl p_s + \int^0_z \nbl\tht(x,y,\xi,t)d\xi + fv^\bot =0.
\end{split}
\end{equation}
\begin{equation}
\label{e:t.n}
 \tht_t + L_2\tht + v\cdot\nbl\tht +w \tht_z =Q ;
\end{equation}
\begin{equation}
\label{e:w}
 w(x, y, z, t) =  \int^0_z \nbl\cdot v(x, y, \xi, t)d\xi,
\end{equation}
\begin{equation}
\label{e:bc.v}
  (v_z+\alf_1 v) |_{z=0} = v |_{z=-h} = v |_{(x,y)\in\prt{D}}=0,
\end{equation}
\begin{equation}
\label{e:bc.w}
  w |_{z=-h} = 0,
\end{equation}
\begin{equation}
\label{e:bc.t}
 (\tht_z+\alf_2 \tht)|_{z=0}= \tht_z|_{z=-h}=\Pn\tht |_{(x,y)\in\prt{D}}=0.
\end{equation}
 The above system of PEs will be solved with suitable initial conditions:
\begin{equation}
\label{e:ic.n}
 v(x, y, z, 0) = v_0(x, y, z), \quad \tht(x, y, z, 0) = \tht_0(x, y, z).
\end{equation}
 Assume $Q$ is independent of time for simplicity of discussion, since results
 to be presented for autonomous case can be extended to non-autonomous case
 with proper modifications.

 The notions of weak and strong solutions were introduced in \cite{l.t.w:92.o},
 where existence of weak solutions was proved, though uniqueness of weak
 solutions is still not resolved yet. Local (in time) existence and uniqueness
 of strong solutions were obtained in \cite{g.m.r:01} and \cite{t.z:04}.
 Global (in time) existence of strong solutions was proved in \cite{c.t:07} for
 the case when $v$ satisfies a set of Neumann type boundary conditions.
 See \cite{k:07} for a different approach.
 For the case when $v$ satisfies the boundary conditions \eqref{e:bc.v}, see
 \cite{k.z:07} for a proof of global regularity of strong solutions. For the
 case when $v$ satisfies a set of related boundary conditions, see \cite{h.k:16}
 for a different approach.
 Uniform boundedness of strong solutions was proved in \cite{j:07} for the case
 with Neumann type boundary conditions and in \cite{k.z:08} for the case with
 physical boundary conditions.  For many more other results and studies of
 3D Primitive Equations, refer \cite{l.t:16},\cite{p.t.z:09} and \cite{t.z:04}.

 An initial motivation of this paper is to study uniform boundedness of $H^2$
 solutions to the 3D primitive equations for the case with physical boundary
 conditions \eqref{e:bc.v}. Global existence of the $H^2$ solutions for the case
 with {\em periodic} boundary condition was proved in \cite{p:06}. The method
 of \cite{p:06} uses extensively integrations by parts. For non-periodic cases,
 the boundary terms arising from those integrations by parts would cause trouble
 for {\em a priori} estimates. The first success overcoming this difficulty was
 achieved in \cite{j.t:15}, which proved global uniform boundedness of $H^2$
 solutions for the 3D Primitive Equations with $v$ satisfying a set of
 {\em Neumann} type boundary conditions. An interesting aspect of this approach
 is that not only it works for both the periodic and the non-periodic case, it
 also requires {\em less demanding} condition on $Q$ than \cite{p:06}. This
 approach was further improved in \cite{j:17} which eliminated all technical
 restrictions {\em completely}.
 Roughly speaking, the main idea of \cite{j:17} and \cite{j.t:15} is obtaining
 uniform boundedness of $\|(v_t,\tht_t)\|_{L^2}$ and $\|v_z\|_{H^1}$ first, then
 using elliptic regularity to get uniform boundedness of $\|(v,\tht)\|_{H^2}$,
 since boundary conditions of $v_t,\tht_t$ and $v_z$ can be used in this
 approach.
 However, it seems the approach of \cite{j:17} and \cite{j.t:15} can not help
 for the case with the physical boundary condition \eqref{e:bc.v}. The reason
 is that uniform boundedness of $\|v_z\|_{H^1}$ is difficult for this case;
 while the nonlinear term $wv_z$ in \eqref{e:v.n} is not controllable with just
 uniform boundedness of $\|(v_t, \tht_t)\|_{L^2}$ alone.
 Another drawback of the method of \cite{j:17} and \cite{j.t:15} is that it
 seems not helpful for uniform estimate of higher regularity, even for the case
 with Neumann boundary conditions.
 
 In this paper, a completely different approach is presented. The new approach
 can prove uniform boundedness of $H^m$ norm of $(v,\tht)$ for any $m\gs 2$, if
 $\prt D$ and $Q$ are sufficiently smooth. Moreover, this approach works not
 only for the case with boundary condition \eqref{e:bc.v}, but also for other
 typical boundary conditions. The main new idea is the use of {\em hydrostatic}
 Leray projector to be defined in Section~\ref{s:helm}. It allows treatment of
 pressure and boundary conditions more convenient. To realize the goal, several
 fundamental results for hydrostatic Leray projector and a new
 {\em anisotrophic} multi-linear estimate in Sobolev space are established.
 With the help of these results, the proofs of the new approach for uniform
 boundedness appear direct and simplified. The results on hydrostatic Leray
 projector are expected to be of general interest and can be helpful for solving
 some other problems about the system of 3D Primitive Equations which appeared
 rather difficult previously, due to non-periodic boundary conditions. For
 example, analyticity in Gevery class for the solutions of 3D Primitive
 Equations with {\em non-periodic} boundary conditions can be established now
 with the help of hydrostatic Leray projector (\cite{j:17.g}).

 The rest of this article is organized as follows:

 In Section~\ref{s:pre}, we give the notations, briefly review the background
 results and recall some important facts which are useful for later analysis.
 We will also state and prove Lemma~\ref{l:ju.2}, which gives a useful new
 anisotrophic multi-linear estimate in Sobolev space and will be used a few
 times in Sections \ref{s:h2} and \ref{s:high}.

 In Section~\ref{s:helm}, we present and prove Theorems \ref{t:h1.g},
 \ref{t:h1.dcp}, \ref{t:P.d} and \ref{t:hso}. These fundamental results on
 hydrostatic Helmholtz decomposition and hydrostatic Leray projector will be
 used in several key steps of the analyses of Section~\ref{s:h2} and
 Section~\ref{s:high}.

 In Section~\ref{s:h2}, the result of global uniform boundedness of $H^2$
 solutions is stated as Theorem~\ref{t:h2} and is proved. Theorem~\ref{t:t} on
 global uniform boundedness of $\|(v_t, \tht_t)\|$ will also be proved.

 In Section~\ref{s:high}, global uniform boundedness for higher regularity
 will be presented as Theorem~\ref{t:hm} and be proved briefly.

\section{Preliminaries}
\label{s:pre}

\noindent
 Notations to be used are basically standard. $C$ and $c$ denote generic
 positive constants, the values of which may vary from one occurrence to
 another. The relation operators $\lsp$ and $\apx$ are used such that, for
 real numbers $A$ and $B$,
\begin{align*}
  A \lsp B \quad & \text{ if and only if }\quad  A\ls C\cdot B,\\
  A \apx B \quad &\text{ if and only if }\quad  c\cdot A\ls B \ls C\cdot A,
\end{align*}
 for some positive constants $C$ and $c$ independent of $A$ and $B$.

 $L^p(\Omg)$ and $L^p(D)$ ($1\ls p\ls\infty$) are the classic Lebesgue $L^p$
 spaces with the norm denoted as
\[ \|\cdot\|_p:=\|\cdot\|_{L^p}, \quad \|\cdot\|:=\|\cdot\|_2.\] 
 $H^m(\Omg)$ and $H^m(D)$ ($m\gs 1$) (with norm $\|\cdot\|_{H^m}$) denote the
 classic Sobolev spaces for $L^2$ functions with weak derivatives up to order
 $m$ also in $L^2$. Notations for vector and scalar function spaces may not be
 distinguished when they are self-evident from the context.

 Define $\ol{\vphi}$ as the vertical average of $\vphi$ on $\Omg$:
\[ \ol{\vphi}(x,y) = \frac{1}{h}\int_{-h}^0 \vphi(x,y,z)dz .\]
 By the H\"older inequality, it is easy to see that, for $\vphi\in L^p(\Omg)$,
\begin{equation}
\label{e:va}
 \|\ol{\vphi}\|_{L^p(\Omg)}  = h^\frac{1}{p}\|\ol{\vphi}\|_{L^p(D)}
 \ls 
 \|\vphi\|_p, \quad \forall p\in[1, +\infty].
\end{equation}
 Define the function spaces:
\begin{equation*}
\begin{split}
 H := H_1\times H_2 := \{ v\in (L^2(\Omg))^2\ | \
  \nbl\cdot\ol{v} =0,\quad \ol{v}\cdot n |_{\prt D}=0 \}\times L^2(\Omg),\\
 V := V_1\times V_2 := \{ v\in (H^1(\Omg))^2\ | \
  \nbl\cdot\ol{v} =0,\quad v |_{\Gmm_b\cap\Gmm_l}=0 \} \times H^1(\Omg).
\end{split}
\end{equation*}
 Notice that, when $\alf_2=0$, $V_2$ is chosen as
  \[ V_2:=\lf\{ \vfi\in H^1(\Omg): \int_\Omg \vfi d\Omg =0 \rt\}. \]
 We use $d\Omg$ and $dD$ to denote $dxdydz$ and $dxdy$ in integrals in $\Omg$
 and $D$ respectively, or we may simply omit them.

 Define the bilinear forms: $a_i: V_i\times V_i \ra \R$, $i=1,2$, such that
\begin{equation*}
\begin{aligned}
 a_1(v, u) &= \int_\Omg \left(\nu_1 \sum_{i=1}^2\nbl v_i\cdot \nbl u_i
 + \mu_1 v_z\cdot u_z \right)d\Omg
 + \alf_1\mu_1\int_{\Gmm_t} v\cdot u\ dD,\\
 a_2(\tht, \eta)& = \int_\Omg \left(\nu_2 \nbl\tht \cdot \nbl\eta 
 + \mu_2 \tht_z \eta_z \right)d\Omg
 + \alf_2\mu_2\int_{\Gmm_t}\tht\eta\ dD,
\end{aligned}
\end{equation*}
 and the linear operators $A_i : V_i \mapsto V_i'$, $i=1,2$, such that
 \begin{equation*}
 \lra{A_1v}{u}=a_1(v, u), \quad \forall v, u \in V_1;
 \quad
 \lra{A_2\tht}{\eta}=a_2(\tht, \eta),
 \quad  \forall \tht, \eta \in V_2,
 \end{equation*}
 where $V_i'$ ($i=1,2$) is the dual space of $V_i$ and $\lra{\cdot}{\cdot}$ is
 the corresponding scalar product between $V_i'$ and $V_i$.
 We also use $\lra{\cdot}{\cdot}$ to denote the inner products in
 $H_1$ and $H_2$.\\
 Define:
 \[ D(A_i) = \{ \phi \in V_i, A_i\phi\in H_i \},\quad i=1,2. \]
 Since $A_i^{-1}$ is a self-adjoint compact operator in $H_i$, by the classic
 spectral theory, the power $A_i^s$ can be defined for any $s\in\R$.
 Then
\[ D(A_i)' = D(A_i^{-1})\]
 is the dual space of $D(A_i)$ and, with the denotation $\Lmd_i=A_i^\hf$,
\[ V_i = D(A_i^\half)=D(\Lmd_i), \quad V_i' = D(A_i^{-\half})=D(\Lmd_i^{-1}).\]
 Moreover,
\[ D(A_i) \subset V_i \subset H_i \subset V_i' \subset D(A_i)', \]
 where the embeddings above are all compact and each space above is dense in the
 one following it.\\
 Define the norm $\|\cdot\|_{V_i}$ by:
 \[ \|\cdot\|_{V_i}^2 = a_i(\cdot,\cdot)
	= \lra{A_i\cdot}{\cdot}
	= \lan A_i^\half\cdot, A_i^\half\cdot \ran
        =\lra{\Lmd_i\cdot}{\Lmd_i\cdot},
	\quad i=1,2. \]
 Then, for any $\phi =(\phi_1, \phi_2)\in V_1$ and $\psi \in V_2$
\begin{equation*}
 \|\phi\| \lsp \|\phi\|_{V_1}, \quad  \|\psi\| \lsp \|\psi\|_{V_2}.
\end{equation*}
 Therefore, for any $\phi =(\phi_1, \phi_2)\in V_1$ and $\psi \in V_2$,
 \begin{equation}
\label{e:nmeqv}
 \|\phi\|_{V_1} \apx \|\phi\|_{H^1}, \quad
 \|\psi\|_{V_2} \apx \|\psi\|_{H^1}.
\end{equation}

 Recall the following definitions of weak and strong solutions:
\begin{defn}
\label{d:soln}
 Suppose $Q \in L^2(\Omg)$, $(v_0,\tht_0)\in H$ and $T>0$.

 The pair $(v,\tht)$ is called a {\em weak solution} of the 3D viscous PEs
 (\ref{e:v.n})-(\ref{e:ic.n}) on the $[0, T]$ if it satisfies
 (\ref{e:v.n})-(\ref{e:ic.n}) in weak sense (see e.g. \cite{l.t.w:92.o} for
 details) and that
\begin{equation*}
 (v, \tht) \in L^\infty([0, T];H) \cap L^2(0, T; V).
\end{equation*}
 Moreover, if $(v_0, \tht_0)\in V$, a weak solution $(v, \tht)$ is called a
 {\em strong solution} of (\ref{e:v.n})-(\ref{e:ic.n}) on the time interval
 $[0, T]$ if, in addition, it satisfies
\begin{equation*}
 (v,\tht) \in L^\infty([0, T]; V) \cap L^2(0, T; D(A_1)\times D(A_2)).
\end{equation*}

\end{defn}

 Recall the following lemma which will be used in the {\em a priori} estimates
 in Sections~\ref{s:h2} and \ref{s:high}. See \cite{c.t:03} and \cite{j:07} for
 the proof.
\begin{lemm}
\label{l:ju.1}
 Suppose that $\nbl v, \vfi\in H^1(\Omg)$ and $\psi\in L^2(\Omg)$. Then,
\begin{equation*}
\lf|\lra{\lf(\int_z^0\nbl\cdot v(x,y,\xi)d\xi\rt) \vfi}{\psi}\rt|
 \lsp \|\nbl v\|^{\hf} \|\nbl v\|_{H^1}^{\hf}
 \|\vfi\|^\hf\|\vfi\|_{H^1}^\hf \|\psi\|.
\end{equation*}
\end{lemm}
 Finally, we prove the following anisotrophic estimate in Sobolev spaces. It
 will be used as well in Sections~\ref{s:h2} and \ref{s:high}.
\begin{lemm}
\label{l:ju.2}
 Suppose $\phi,\psi\in H^1(\Omg)$ and $\vfi\in L^2(\Omg)$. Then,
\begin{equation*}
\begin{split}
 \int_\Omg |\phi\psi\vfi| \lsp&
 \|\phi\|^\qt(\|\phi\|+\|\phi_z\|)^\qt(\|\phi\|+\|\nbl\phi\|)^\hf\\
 &\times \|\psi\|^\qt(\|\psi\|+\|\psi_z\|)^\qt(\|\psi\|+\|\nbl\psi\|)^\hf
 \|\vfi\|.
\end{split}
\end{equation*}
\end{lemm}

\begin{prf}
\begin{equation*}
\begin{split}
 \int_\Omg |\phi\psi\vfi|
\ls & \int_{-h}^0 \|\phi\|_{L^4(D)}\|\psi\|_{L^4(D)}\|\vfi\|_{L^2(D)}\ dz\\
\lsp& \int_{-h}^0 \|\phi\|_{L^2(D)}^\hf
 \lf(\|\phi\|_{L^2(D)}+\|\nbl\phi\|_{L^2(D)}\rt)^\hf\\
 &\quad \x \|\psi\|_{L^2(D)}^\hf
 \lf(\|\psi\|_{L^2(D)}+\|\nbl\psi\|_{L^2(D)}\rt)^\hf \|\vfi\|_{L^2(D)}\ dz \\
\ls& \|\phi\|_{L_z^\infty(L_D^2)}^\hf\|\psi\|_{L_z^\infty(L_D^2)}^\hf \\
 & \x \int_{-h}^0 \lf(\|\phi\|_{L^2(D)}+\|\nbl\phi\|_{L^2(D)}\rt)^\hf\\
 & \qquad\x
  \lf(\|\psi\|_{L^2(D)}+\|\nbl\psi\|_{L^2(D)}\rt)^\hf\|\vfi\|_{L^2(D)}\ dz\\
\lsp& \|\phi\|_{L_D^2(L_z^\infty)}^\hf\|\psi\|_{L_D^2(L_z^\infty)}^\hf
  (\|\phi\|+\|\nbl\phi\|)^\hf(\|\psi\|+\|\nbl\psi\|)^\hf\|\vfi\|
\end{split}
\end{equation*}
 Notice that
\begin{equation*}
\begin{split}
 \|\phi\|_{L_D^2(L_z^\infty)}^2 =& \int_D \|\phi\|_{L_z^\infty}^2\ dD\\
 \lsp& \int_D \|\phi\|_{L_z^2}(\|\phi\|_{L_z^2}+\|\phi_z\|_{L_z^2})\ dD
 \lsp \|\phi\|_2(\|\phi\|_2+\|\phi_z\|_2).
\end{split}
\end{equation*}
 This finishes the proof of Lemma~\ref{l:ju.2}.
\end{prf}

\section{Hydrostatic Decomposition}
\label{s:helm}

 First, we prove the following Helmholtz type hydrostatic decomposition.
 Recall that the space $H_1$ was defined in Section~\ref{s:pre}.
\begin{thrm}
\label{t:h1.g}
 \[ (L^2(\Omg))^2 =  H_1 \oplus G, \]
 where
\[ G :=\{ u\in (L^2(\Omg))^2 \ |\  u_z=0,\ u = \nbl q,\quad q\in H^1(D) \}. \]
\end{thrm}
\begin{prf}
 Let $(H_1)^\bot$ be the orthogonal completement of $H_1$ with respect
 to the inner product in $(L^2(\Omg))^2$, i.e.
\[ (L^2(\Omg))^2 = H_1 \oplus (H_1)^\bot. \]
 We need only to show that $G=(H_1)^\bot$.

 {\bf Claim 1}: $G \subset (H_1)^\bot$.\\
 {\em Proof of Claim 1}:
 Let $u\in G $. Thus, $u=\nbl q$ for some $q\in H^1(D)$ and $u_z=0$.
 Let $v \in {\cal V}_1$, where
\[ {\cal V}_1
 := \{ v\in ({\cal C}^\infty(\Omg))^2\ |\
  v \text{ vanishes in a neighborhood of } \Gmm_b\cup\Gmm_l,\
  \nbl\cdot \bar{v}=0 \}. \]
 Then,
\[ \lra{u}{v} = h \int_D (\nbl q)\cdot \bar{v} dD
  = - h \int_D q \nbl\cdot\bar{v} = 0. \]
  Recall that ${\cal V}_1$ is dense in $V_1$ (see \cite{l.t.w:92.a}). Hence,
 ${\cal V}_1$ is dense in $H_1$, since $V_1$ is dense in $H_1$.
 Therefore, $u\in (H_1)^\bot$. This proves Claim 1.

 {\bf Claim 2}: $(H_1)^\bot \subset G$.\\
 {\em Proof of Claim 2}:
 Let $u=(u_1, u_2)\in (H_1)^\bot$. Then, $u\in (L^2(\Omg))^2$ and
\[ \lra{u}{v} = 0, \quad \forall v\in H_1.\]
 {\em Step 1}.
 Choose the special $v =(\vfi_z, 0)$ with $\vfi\in {\cal C}_0^\infty(\Omg)$.
 Then, $v \in H_1$, since $\bar{v}=(0,0)$. Thus,
\[ \lra{u_1}{\vfi_z} = \lra{u}{v}=0, \quad
  \forall \vfi\in {\cal C}_0^\infty(\Omg). \]
 Therefore, $\Pz u_1=0$ as the weak derivative of $u_1$. Similarly, $\Pz u_2=0$.
 Thus,
\[ u_z=0, \quad \text{and}\ u\in (L^2(D))^2. \]
 {\em Step 2}.
 Choose $v\in (L^2(\Omg))^2$ such that  $v_z=0$ and
\[ \nbl\cdot v =0, \quad n\cdot v =0|_{\prt D}=0. \]
 Since $v\in H_1$, we have $\lra{u}{v}=0$ for all such $v$'s. Thus, there exists
 a $q\in H^1(D)$, such that
\[ u = \nbl q. \]
 Therefore, $u\in G$. This proves Claim 2.
\end{prf}
 
 An immediate application of Theorem~\ref{t:h1.g} is the definition of the
 {\em hydrostatic Leray projector} $\Pf$ as the orthogonal projection of
 $(L^2(\Omg))^2$ onto $H_1$ with repect to the inner product of $(L^2(\Omg))^2$.

 Next, we give the following refined decomposition result for $(L^2(\Omg))^2$:
\begin{thrm}
\label{t:h1.dcp}
 Let $D$ be an open bounded set of class ${\cal C}^2$. Then,
\begin{equation}
\label{e:dec.s}
 (L^2(\Omg))^2 = H_1\oplus G_1\oplus G_2,
\end{equation}
 where $H_1$, $G_1$ and $G_2$ are mutually orthogonal spaces and
\begin{align*}
 G_1
 :=& \{ u\in (L^2(\Omg))^2\ |\ u_z=0,\ u=\nbl q,\ q\in H^1(D),\ \Dlt q =0 \},\\
 G_2 :=& \{ u\in (L^2(\Omg))^2\ |\ u_z=0,\ u=\nbl q,\ q\in H_0^1(D) \}.
\end{align*}
 Moreover, the following decomposition is valid:
 \begin{equation}
\label{e:dec.u}
 u= \Pf u + \nbl(q_1+q_2), \quad \forall u\in (L^2(\Omg)),
\end{equation}
 where
\begin{align}
\label{e:q1}
 \Dlt q_1 = 0, \quad & \Pn q_1|_{\prt D} = n\cdot(\bar{u}-\nbl q_2), \\
\label{e:q2}
 \Dlt q_2 = \nbl\cdot\bar{u}\in H^{-1}(D), \quad & q_2\in H_0^1(D).
\end{align}
 \end{thrm}
\begin{prf}
 In the following, we prove \eqref{e:dec.u}, from which \eqref{e:dec.s}
 follows as well. For any given $u\in (L^2(\Omg))^2$, the Dirichlet problem
 \eqref{e:q2} has a unique solution $q_2$. Therefore, the Neumann problem
 \eqref{e:q1} is well defined. Moreoever, by the Stokes formula and
 \eqref{e:q2}, we have
\[ \int_{\prt D} n\cdot(\bar{u}-\nbl q_2)\ ds
 = \int_D \nbl\cdot({\bar{u}-\nbl q_2})\ dD
 = \int_D (\nbl\bar{u} -\Dlt q_2)\ dD =0. \]
 Then, the Neumann problem \eqref{e:q1} has a solution $q_1\in H^1(D)$, which is
 unique up to an additive constant. Obviously,
\[ \nbl q_i \in Q_i,\ i=1, 2. \]
 It is easy to see that $\nbl q_1$ and $\nbl q_2$ are orthogonal in
 $(L^2(\Omg))^2$, since
\begin{equation*}
\begin{split}
\int_\Omg \nbl q_1\cdot \nbl q_2\ d\Omg
 =& h \int_D \nbl q_1\cdot \nbl q_2\ dD\\
 =& h \int_{\prt D} (\Pn q_1) q_2 ds - h \int_D (\Dlt q_1) q_2\ dD
 =0,
\end{split}
\end{equation*}
 where the Stokes formula and the definitions of $q_1$ and $q_2$ given by
 \eqref{e:q1} and \eqref{e:q2} are used. Finally, we show that
\[ u_* : = u- \nbl(q_1+q_2) \in H_1. \]
 Notice that, by \eqref{e:q1} and \eqref{e:q2},
\[ \nbl\cdot \ol{u_*} = \nbl\cdot\bar{u} - \Dlt q_1 -\Dlt q_2
 = \nbl \cdot \bar{u} -\Dlt q_2=0, \]
 and
\[ \ol{u_*}\cdot n|_{\prt D}= \lf(n\cdot \bar{u} + \Pn q_1 +\Pn q_2\rt)|_{\prt D}
 =0. \]
 Thus, $u_*\in H_1$. This proves \eqref{e:dec.u} and hence \eqref{e:dec.s}
 as well.
\end{prf}

 The importance of \eqref{e:dec.u} is that it gives more concrete information
 about the decomposition \eqref{e:dec.s}. Especially, it provides detailed
 relation between $\Pf u$ and $u$, via \eqref{e:q1} and \eqref{e:q2}.
 An immediate and very important consequence of \eqref{e:dec.u} is the
 following result.
\begin{thrm}
\label{t:P.d} For any $u\in (L^2(\Omg))^2$, the following statements are valid:
 \begin{enumerate}
\item[(a)] Suppose $|u_z|\in L_{loc}^1(\Omg)$. Then, in the sense of
 distribution,
\begin{equation*}
 (\Pf u)_z = u_z.
\end{equation*}
 \item[(b)] Suppose $|u|, |\nbl u|\in L^r(\Omg)$, $r\in(1,\infty)$ and
 $\prt D\in {\mathcal C}^2$. Then, there exists a constant $c=c(r,D)>0$
 such that:
\begin{equation*}
 \||\nbl(\Pf u)|\|_{L^r(\Omg)}
 \ls c(r,D) (\|u\|_{(L^r(\Omg))^2}+\||\nbl u|\|_{L^r(\Omg)}).
\end{equation*}
 \item[(c)] Suppose $u\in (W^{1,r}(\Omg))^2$, $r\in(1,\infty)$ and
 $\prt D\in {\mathcal C}^2$. Then, there exists a constant $c=c(r,D)>0$
 such that:
\[ \||(\nbl, \Pz)\Pf u|\|_{L^r(\Omg)} 
 \ls c(r,D) \|u\|_{(W^{1,r}(\Omg))^2}. \]
\item[(d)] Suppose $u\in (W^{m,r}(\Omg))^2$, $r\in(1,\infty)$ and
 $\prt D\in {\mathcal C}^{m+1}$. Then, there exists a constant $c=c(r,D)>0$
 such that:
\[ \||\nbl_3^m\Pf u|\|_{L^r(\Omg)} 
 \ls c(r,D) \|u\|_{(W^{m,r}(\Omg))^2}. \]
\end{enumerate}
\end{thrm}
\begin{prf}
 Theorem~\ref{t:P.d} (a) follows from \eqref{e:dec.u} immediately.
 By \eqref{e:dec.u}, we also have
\[ \||\nbl(\Pf u)|\|_{L^r(\Omg)} \ls\|\nbl u\|_{(L^r(\Omg))^2}
 + h^\frac{1}{r}
 \lf(\||\nbl^2 q_1|\|_{L^r(D)} + \||\nbl^2 q_2|\|_{L^r(D)}\rt). \]
 From \eqref{e:q1}, it follows that
\begin{equation*}
\begin{split}
 \||\nbl^2 q_1|\|_{L^r(D)}
\lsp& \||\bar{u}-\nbl q_2|\|_{W^{1-1/r,r}(\prt D)}\\
\lsp& \|\bar{u}\|_{(W^{1,r}(D))^2} +\|\nbl q_2\|_{(W^{1,r}(D))^2}\\
\lsp& h^{-\frac{1}{r}}(\|u\|_{(L^r(\Omg))^2}+\||\nbl u|\|_{L^r(\Omg)})
 +\|\nbl q_2\|_{(W^{1,r}(D))^2}
\end{split}
\end{equation*}
 From \eqref{e:q2}, it follows that
\[ \|\nbl q_2\|_{(W^{1,r}(D))^2}
 \lsp \|\nbl\cdot \bar{u}\|_{L^r(D)}
 \lsp h^{-\frac{1}{r}}\||\nbl u|\|_{L^r(\Omg)}. \]
 In the above estimates, \eqref{e:va} has been used. Therefore,
\[ \||\nbl(\Pf u)|\|_{(L^r(\Omg))^2}
 \lsp \|u\|_{L^r(\Omg)}+\||\nbl u|\|_{L^r(\Omg)}. \]
 This proves Theorem~\ref{t:P.d} (b). Theorem~\ref{t:P.d} (c) is an immediate
 consequence of (a) and (b). Proof of Theorem~\ref{t:P.d} (d) is similar to that
 of (c).
\end{prf}

 The last important result of the section is about the relation between the
 anisotrophic Laplace operator $L_1$ and the {\em hydrostatic Stokes operator}
 $A_1$ through the hydrostatic Leray projector $\Pf$.  Recall that the operator
 $A_1$ was already defined in Section~\ref{s:pre}.
\begin{thrm}
\label{t:hso}
 As an isomorphism from $V_1$ onto $(V_1)'$, $A_1$ satisfies:
 \[ A_1u = \Pf L_1 u, \quad \forall u\in V_1. \]
\end{thrm}

\begin{prf} It suffices to prove $A_1 = \Pf L_1$ only. Suppose $u\in D(A_1)$
 and $\vfi\in V_1$. Then, $L_1u\in (L^2(\Omg))^2$. By Theorem~\ref{t:h1.g},
 there exists a $q\in G$ such that
\[ \lra{\Pf L_1 u}{\vfi} =\lra{L_1 u -\nbl q}{\vfi} = \lra{L_1 u}{\vfi}
 = -\int_\Omg \lf[\nu_1\Dlt u\cdot\vfi +\mu_1u_{zz}\cdot\vfi\rt] d\Omg.\]
 Since $u\in D(A_1)$ and $\vfi\in V_1$, we have
\begin{align*}
 \lra{\Dlt u}{\vfi} =& \int_{-h}^0 \lf(\int_D\Dlt u\cdot\vfi dD\rt)dz
 = \int_{-h}^0 \lf(\int_D\sum_{i=1}^2\Dlt u_i\vfi_i dD\rt)dz\\
 =& \int_{-h}^0 \lf( \int_{\prt D}\sum_{i=1}^2  \prt_n u_i\cdot\vfi_i ds
  - \int_D\sum_{i=1}^2\nbl u_i\cdot\nbl \vfi_i dD \rt) dz\\
 =&-\int_\Omg\sum_{i=1}^2\nbl u_i\cdot\nbl \vfi_i d\Omg,\\
 \lra{u_{zz}}{\vfi} =& \int_D \lf(\int_{-h}^0u_{zz}\cdot\vfi dz\rt)dD
 = \int_D \lf( u_z\cdot\vfi\Big|_{z=-h}^0
  - \int_{-h}^0 u_z\cdot\vfi_z dz \rt) dD\\
 =&-\alf_1\int_{\Gmm_t} u\cdot\vfi\ dD-\int_\Omg u_z\cdot\nbl \vfi_z d\Omg.
\end{align*}
 Thus,
\[ \lra{\Pf L_1 u}{\vfi} =\lra{A_1 u}{\vfi}, \quad \forall\ \vfi\in V_1, \]
 that is
\[ \Pf L_1 u = A_1 u, \quad \forall\ u\in D(A_1). \]
 Notice that $D(A_1)$ is dense is $V_1$. Therefore, considered as an element
 of $(V_1)'$,
\[ \Pf L_1 u = A_1 u, \quad \forall\ u\in V_1. \]
\end{prf}

\section{Uniform Boundedness of $\|(A_1v,A_2\tht)\|$}
\label{s:h2}
 In this section, we prove {\em uniform} boundedness of $\|(A_1v,A_2\tht)\|$
 and existence of a bounded absorbing set of it in $\R_+$. More precisely, we
 prove the following theorem.
\begin{thrm}
\label{t:h2}
 Suppose $Q\in L^2(\Omg)$ and  $(v_0,\tht_0)\in D(A_1)\x D(A_2)$. Let $(v,\tht)$
 be the unique strong sloution of problem \eqref{e:v.n}-\eqref{e:ic.n}. Then,
\[ (A_1v, A_2\tht)\in L^\infty(0, \infty; H),\ 
  (A_1v, \tht_t) \in L^2(0, T; V),\ \forall\ T\in(0, \infty).\]
 Moreover, there exists a bounded absorbing set of $\|(A_1v,A_2\tht)\|$ in
 $R_+$.
\end{thrm}
\begin{prf}
 {\em Step 1}. Estimate of $\|A_1v\|$.\\
 First, apply the hydrostatic Leray projector $\Pf$ to \eqref{e:v.n} and use
 Theorem~\ref{t:hso} to get
\begin{equation}
\label{e:v.p}
 v_t+A_1v= -\Pf\lf[ (v\cdot\nbl)v+wv_z+\int_z^0 \nbl\tht\ d\xi +fv^\bot \rt].
 \end{equation}
 Next, apply $A_1$ to \eqref{e:v.p} and then take inner product with $A_1v$ to
 arrive at
\begin{equation*}
\begin{split}
\hf \dt{\|A_1v\|^2}&+\|A_1^{\frac{3}{2}}v\|^2\\
=&-\lra{A_1\Pf\lf[(v\cdot\nbl)v+wv_z+\int_z^0\nbl\tht\ d\xi+fv^\bot\rt]}{A_1v}\\
\ls&\lf\|A_1^\hf\Pf\lf[(v\cdot\nbl)v+wv_z+\int_z^0\nbl\tht\ d\xi+fv^\bot\rt]\rt\|
 \|A_1^{\frac{3}{2}}v\|.
\end{split}
\end{equation*}
 Notice that $\Pf$ does {\em not} commute with $A_1$ due to no-slip boundary
 condition. Therefore, Theorem~\ref{t:P.d} is crucial here in dealing with the
 right-hand side of the above inequality, from which it then follows that
\begin{equation*}
\begin{split}
\hf \dt{\|A_1v\|^2} &+\|A_1^{\frac{3}{2}}v\|^2\\
\lsp& \lf\|P\lf[(v\cdot\nbl)v+wv_z+\int_z^0\nbl\tht\ d\xi+fv^\bot \rt]\rt\|_{H^1}
      \|A_1^{\frac{3}{2}}v\|\\
\lsp& \lf\|(v\cdot\nbl)v+wv_z+\int_z^0\nbl\tht\ d\xi+fv^\bot\rt\|_{H^1}
\|A_1^{\frac{3}{2}}v\|,
\end{split}
\end{equation*}
 where we have also used norm equivalence \eqref{e:nmeqv}. Noticing that
 $v|_{\Gmm_b\cup\Gmm_l}=w|_{\Gmm_t\cup\Gmm_b}=0$, we have
\begin{equation*}
\begin{split}
 &\lf\|(v\cdot\nbl)v+wv_z+\int_z^0\nbl\tht\ d\xi+fv^\bot\rt\|_{H^1}\\
\lsp& \| \nbl_3 [(v\cdot\nbl)v+wv_z] \|
  + \lf\| \int_z^0\nbl\tht\ d\xi \rt\|_{H^1} + \|v\|_{V_1}\\
\lsp& \| \nbl_3 [(v\cdot\nbl)v+wv_z]\| + \|v\|_{V_1} + \|A_2\tht\|.
\end{split}
\end{equation*}
 Therefore,
\begin{equation}
\label{e:Av}
 \dt{\|A_1v\|^2}+\|A_1^{\frac{3}{2}}v\|^2
 \lsp \| \nbl_3[(v\cdot\nbl)v+wv_z ] \|^2 + \|v\|_{V_1}^2+\|A_2\tht\|^2.
\end{equation}
 Now, we estimate the first term on the right side of \eqref{e:Av}.
 By Lemma~\ref{l:ju.2}, we have
\begin{align*}
 \|\nbl[(v\cdot\nbl)v]\|
\ls& \||v||\nbl^2 v|\| + \||\nbl v|^2\|\\
\lsp& \|v\|^\qt\|v\|_{V_1}^{\frac{3}{4}}
 \|\nbl^2 v\|^\qt\|\nbl^2v\|_{V_1}^{\frac{3}{4}}
 + \|\nbl v\|^\hf\|\nbl v\|_{V_1}^{\frac{3}{2}}\\
\lsp & \|v\|^\qt \|v\|_{V_1}^{\frac{3}{4}}
 \|A_1v\|^\qt\|A_1^{\frac{3}{2}}v\|^{\frac{3}{4}}
 + \|v\|_{V_1}^\hf\|A_1v\|^{\frac{3}{2}},\\
 \|\Pz[(v\cdot\nbl)v]\|
\lsp& \|(v\cdot \nbl) v_z \| + \|(v_z\cdot\nbl) v \|\\
\lsp& \|v\|^\qt\|v\|_{V_1}^{\frac{3}{4}}
 \|\nbl v_z\|^\qt\|\nbl v_z\|_{V_1}^{\frac{3}{4}}
 + \|v_z\|^\qt \|v_z\|_{V_1}^{\frac{3}{4}}
 \|\nbl v\|^\qt \|\nbl v\|_{V_1}^{\frac{3}{4}}\\
\lsp& \|v\|^\qt\|v\|_{V_1}^{\frac{3}{4}}
 \|A_1v\|^\qt\|A_1^{\frac{3}{2}}v\|^{\frac{3}{4}}
 +\|v\|_{V_1}^\hf\|A_1v\|^{\frac{3}{2}}.
\end{align*}
 Therefore,
\begin{equation}
\label{e:v.nv.h1}
\begin{split}
 \|\nbl_3[(v\cdot\nbl)v]\|^2
\lsp & \|v\|^\hf \|v\|_{V_1}^{\frac{3}{2}}
 \|A_1v\|^\hf\|A_1^{\frac{3}{2}}v\|^{\frac{3}{2}} + \|v\|_{V_1}\|A_1v\|^3\\
\lsp & C_\veps \|v\|^2 \|v\|_{V_1}^6\|A_1v\|^2+ \|v\|_{V_1}\|A_1v\|^3
 +\veps \|A_1^{\frac{3}{2}}v\|^2.
\end{split}
\end{equation}
 Next, we use \eqref{e:w} and Lemma~\ref{l:ju.1} and to get
\begin{equation*}
\begin{split}
    \|\nbl(w v_z)\|
\ls& \|w (\nbl v_z)\| +\|(\nbl w)\ox v_z\|\\
\lsp& \|\nbl v\|^\hf\|\nbl v\|_{V_1}^\hf \|\nbl v_z\|^\hf\|\nbl v_z\|_{V_1}^\hf\\
 &+ \|\nbl^2 v\|^\hf\|\nbl^2 v\|_{V_1}^\hf \|v_z\|^\hf\|v_z\|_{V_1}^\hf\\
\lsp & \|v\|_{V_1}^\hf\|A_1v\|\|A_1^{\frac{3}{2}}v\|^\hf.
\end{split}
\end{equation*}
 Similar to the above estimates, we use Lemma~\ref{l:ju.1} and
 Lemma~\ref{l:ju.2} to get
\begin{equation*}
\begin{split}
 \|\Pz(w v_z)\|
\ls& \|w v_{zz}\| +\| (\nbl\cdot v) v_z\|\\
\lsp& \|\nbl v\|^\hf\|\nbl v\|_{V_1}^\hf \|v_{zz}\|^\hf\|v_{zz}\|_{V_1}^\hf\\
 &+ \|\nbl v\|^\qt\|\nbl v\|_{V_1}^{\frac{3}{4}}
 \|v_z\|^\qt\|v_z\|_{V_1}^{\frac{3}{4}}\\
\lsp& \|v\|_{V_1}^\hf\|A_1v\|\|A_1^{\frac{3}{2}}\|^\hf
 + \|v\|_{V_1}^\hf\|A_1v\|^{\frac{3}{2}}.
\end{split}
\end{equation*}
 Therefore,
\begin{equation}
\label{e:v.wz.h1}
\begin{split}
 \|\nbl_3(wv_z)\|^2
\lsp & \|v\|_{V_1}\|A_1 v\|^2\|A_1^{\frac{3}{2}}v\|
  + \|v\|_{V_1}\|A_1v\|^3\\
\lsp& C_\veps\|v\|_{V_1}^2\|A_1 v\|^4
  + \|v\|_{V_1}\|A_1v\|^3 +\veps \|A_1^{\frac{3}{2}}v\|^2.
\end{split}
\end{equation}
 Combining \eqref{e:Av}, \eqref{e:v.nv.h1} and \eqref{e:v.wz.h1} with $\veps>0$
 sufficiently small, we have
\begin{equation}
\label{e:Av.2}
\begin{split}
\dt{\|A_1v\|^2}&+\|A_1^{\frac{3}{2}}v\|^2\\
 \lsp& \|v\|^2 \|v\|_{V_1}^6\|A_1v\|^2+ \|v\|_{V_1}\|A_1v\|^3\\
  &+ \|v\|_{V_1}^4\|A_1 v\|^4 + \|v\|_{V_1}^2+ \|A_2\tht\|^2\\
 \lsp& (\|v\|^2 \|v\|_{V_1}^6 + \|v\|_{V_1}\|A_1v\| +\|v\|_{V_1}^4\|A_1 v\|^2)
 \|A_1v\|^2\\
 &+  \|v\|_{V_1}^2+ \|A_2\tht\|^2.
\end{split}
\end{equation}
 Now, an application of Gronwall lemma to \eqref{e:Av.2} immediately yields
 the following boundedness:
\[ A_1v\in L^\infty(0,T; H_1) \cap L^2(0,T; V_1),\ \forall\ T\in(0,\infty). \]
 Let
\[ y(t):= \|A_1 v(t)\|^2, \quad h(t):= \|v\|_{V_1}^2+ \|A_2\tht\|^2. \]
\[ g_1(t):= \|v\|^2 \|v\|_{V_1}^6
 + \|v\|_{V_1}\|A_1v\| +\|v\|_{V_1}^4\|A_1 v\|^2.\]
 The following integrals with $t\gs0$
\[ \int_t^{t+1}y(s) ds,\ \int_t^{t+1}h(s) ds,\ \int_t^{t+1}g_1(s) ds \]
 are all uniformly bounded for strong solutions. Then, an application of the
 uniform Gronwall lemma (see \cite{f.p:67} and \cite{t:inf}) to \eqref{e:Av.2}
 yields the folloing uniform boundedness:
\[ A_1v\in L^\infty(0,\infty; H_1)\cap L^2(0,T;V_1),\ \forall\ T\in(0,\infty). \]
 Moreover, there is a bounded absorbing set for $A_1v$ in $H_1$. This finishes
 proof of uniform boundedness of $\|A_1v\|$.

 {\em Step 2}. Estimate of $\|A_2\tht\|$.

 We will estimate $\|A_2\tht\|$ in a way different from that of {\em Step 1}.
 This is because the method of {\em Step 1} will not work for $Q\in L^2(\Omg)$.

 Take inner product of \eqref{e:t.n} with $A_2\tht_t$, we get
\begin{equation}
\begin{split}
\label{e:At}
 \|A_2^\hf\tht_t\|^2 + \hf \dt{\|A_2\tht\|^2}
=& \lra{Q}{A_2\tht_t} + \lra{A_2^\hf[ v\cdot\nbl\tht +w\tht_z]}{A_2^\hf\tht_t}\\
=& \dt{\lra{Q}{A_2\tht}}+\lra{A_2^\hf[ v\cdot\nbl\tht +w\tht_z]}{A_2^\hf\tht_t}.
\end{split}
\end{equation} 
 Therefore,
\begin{equation*}
 \|A_2^\hf\tht_t\|^2 + \dt{\|A_2\tht\|^2}
\ls 2\dt{\lra{Q}{A_2\tht}} +\|A_2^\hf[ v\cdot\nbl\tht +w\tht_z]\|^2.
\end{equation*}
 Noticing that $v|_{\Gmm_b\cup\Gmm_l}=w|_{\Gmm_t\cup\Gmm_b}=0$, we have
\begin{equation}
\label{e:At.1}
 \|A_2^\hf\tht_t\|^2 + \dt{\lf[\|A_2\tht\|^2-2\lra{Q}{A_2\tht}\rt]}
 \lsp \|\nbl_3[ v\cdot\nbl\tht +w\tht_z]\|^2.
\end{equation}
 Now, we estimate the right-hand side of \eqref{e:At.1}. We have
\begin{align*}
 \|\nbl_3(v\cdot\nbl\tht)\|
\lsp& \|(v\cdot\nbl)\nbl_3\tht \| + \||\nbl_3 v||\nbl\tht|\|
\lsp \|A_1v\| \|A_2\tht\| + \|\nbl_3 v\|_3\|\nbl\tht\|_6,\\
 \|\nbl(w\tht_z)\|
\lsp& \|w\nbl\tht_z \| + \|\tht_z\nbl w \|
\lsp \|w\|_\infty \|A_2\tht\| + \|\nbl^2v\|_3\|\tht_z\|_6,\\
 \|\Pz(w\tht_z)\|
\lsp& \|w\tht_{zz}\| + \|\tht_z|\nbl v| \|
\lsp \|w\|_\infty \|A_2\tht\| + \|\nbl v\|_3\|\tht_z\|_6.
\end{align*}
 By interpolation inequality and Sobolev embedding, we have
\[ \|\nbl_3 v\|_3 \lsp \|v\|_{V_1}^\hf\|A_1v\|^\hf, \
 \|\nbl^2 v\|_3 \lsp \|A_1v\|^\hf\|A_1^{\frac{3}{2}}v\|^\hf, \
 \|\nbl\tht\|_6,\|\tht_z\|_6 \lsp \|A_2\tht\|.\]
 Using \eqref{e:w}, we have
\[ |w(x,y,z)| = \lf|\int_{-h}^z \nbl\cdot v(x,y,\xi)\ d\xi\rt|
 \ls \int_{-h}^0 |\nbl v(x,y,\xi)|\ d\xi. \]
 By Agmon's inequality,
\begin{equation*}
\begin{split}
 |w(x,y,z)| \ls&
 \int_{-h}^0 \|\nbl v(\cdot,\cdot,\xi)\|_{(L^2(D))^2}^\hf
 \|\nbl v(\cdot,\cdot,\xi)\|_{(H^2(D))^2}^\hf\ d\xi\\
 \lsp& h^\hf\|v\|_{V_1}^\hf \|A_1^{\frac{3}{2}}v\|^\hf.
\end{split}
\end{equation*}
 Thus,
\[ \|w\|_\infty \ls h^\hf\|v\|_{V_1}^\hf \|A_1^{\frac{3}{2}}v\|^\hf. \]
 Collecting all the above estimates after \eqref{e:At.1}, we obtain
\[  \|\nbl_3(v\cdot\nbl\tht + w\tht_z)\|^2 \lsp g_2(t)\|A_2\tht\|^2,\]
 where
\[ g_2(t):= \|v\|_{V_1}\|A_1v\|+\|A_1v\|^2+ h\|v\|_{V_1}\|A_1^{\frac{3}{2}}v\|
 +\|A_1v\| \|A_1^{\frac{3}{2}}v\|. \]
 Then, \eqref{e:At.1} implies
\begin{equation}
\label{e:At.2}
 \|A_2^\hf\tht_t\|^2 + \dt{\|A_2\tht-Q\|^2}
 \lsp  g_2(t) (\|A_2\tht-Q\|^2 +\|Q\|^2).
\end{equation}
 By \eqref{e:Av.2} and uniform boundedness of $\|A_1v\|$, we have
\[ g_2 \in L^2(0, T), \quad \forall\ T\in(0,\infty). \]
 Thus, we conclude immediately from \eqref{e:At.2} by Gronwall lemma that
 \[ A_2\tht \in L^\infty(0, T; H_2),\ \tht_t\in L^2(0, T; V_2),\
 \forall\ T\in(0,\infty). \]
 Moreover, by \eqref{e:Av.2}, we see that $\int_t^{t+1} g_2(s)\ ds$
 is uniformly bounded for all $t\gs 0$ and it is also obvious that
\[ \int_t^{t+1} \|A_2\tht(s)-Q\|^2\ ds
  \ls \|Q\|^2 + \int_t^{t+1}\|A_2\tht(s)\|^2\ ds \]
 is uniformly bounded for all $t\gs 0$, since $(v,\tht)$ is a strong solution.
 Thus, an application of uniform Gronwall lemma to \eqref{e:At.2} yields an
 absorbing set for $\|A_2\tht\|$ in $\R_+$ and the uniform boundedness:
\[ A_2\tht \in L^\infty(0, \infty; H_2). \]
 \end{prf}
 
\begin{rmk}
 If $Q\in H^1(\Omg)$, then, it is easy to see that one can use the method in
 {\em Step 1} of the proof of Theorem~\ref{t:h2} to estimate $\|A_2\tht\|$ and
 obtain the following additional regularity for $\tht$:
\[ A_2\tht \in L^2(0, T; V_2),\ \forall\ T\in(0, \infty). \]
\end{rmk}
 
 With Theorem~\ref{t:h2}, we obtain the following Theorem~\ref{t:t}. Notice
 that only $(u_0,\tht_0)\in D(A)$ is assumed here. Especially,
 $(u_t(0),\tht_t(0))\in H$ is {\em not} assumed, {\em nor} is it assumed that
 the system of 3D primitive equations is valid at $t=0$, as the later implies
 $(u_t(0),\tht_t(0))\in H$ when $(u_0,\tht_0)\in D(A)$.

\begin{thrm}
\label{t:t}
 Suppose $Q\in L^2(\Omg)$ and  $(v_0,\tht_0)\in D(A_1)\x D(A_2)$. Let $(v,\tht)$
 be the unique strong solution of the problem \eqref{e:v.n}-\eqref{e:ic.n}.
 Then,
 \[ (v_t,\tht_t)\in L^\infty(0, \infty; H)\cap L^2(0, T; V),\
  \forall\ T\in(0, \infty). \]
 Moreover, there exists a bounded absorbing set of $\|(v_t,\tht_t)\|$ in $\R_+$.
\end{thrm}

\begin{prf}

 By \eqref{e:v.p} and Lemma~\ref{l:ju.1}, we have
\begin{equation}
\label{e:vt}
\begin{split}
 \|v_t\| \lsp& \|A_1v\| + \|(v\cdot\nbl)v\|+\|wv_z\| +\|\nbl\tht\|+\|v\|\\
\lsp& \|A_1v\|+ \|v\|_6\|\nbl v\|_3
 +\|\nbl v\|^\hf\|\nbl v\|_{V_1}^\hf \|v_z\|^\hf\|v_z\|_{V_1}^\hf
 +\|\nbl\tht\|+\|v\|\\
\lsp&(1+\|v\|_{V_1}) \|A_1v\|+\|v\|_{V_1}^{\frac{3}{2}}\|A_1v\|^\hf
  +\|\nbl\tht\|+\|v\|.
\end{split}
\end{equation}
 By \eqref{e:t.n}, we have
\begin{equation}
\label{e:tt}
\begin{split}
 \|\tht_t\| \ls& \|A_2\tht\|+ \|(v\cdot\nbl)\tht\| + \|w\tht_z\|\\
\lsp& \|A_2\tht\| + \|A_1v\|\|\nbl\tht\|+\|A_1v\|\|A_2\tht\|.
\end{split}
\end{equation}
 Thus, uniform boundedness of $\|(v_t,\tht_t)\|$ and existence of an absorbing
 ball of it in $\R_+$ follow from \eqref{e:vt}, \eqref{e:tt} and
 Theorem~\ref{t:h2}.

 Denote
\[ u:=v_t=\prt_t v, \quad \zeta:=\tht_t=\prt_t\tht. \]
 By \eqref{e:v.n} and \eqref{e:t.n}, we have
\begin{equation}
\label{e:vt.2}
\begin{split}
 u_t + L_1u &+ (u\cdot\nbl)v +(v\cdot\nbl)u  + w_tv_z + w u_z\\
& + \nbl (p_s)_t +\int_z^0\nbl\zeta(x,y,\xi,t)d\xi + fu^\bot =0,
\end{split}
\end{equation}
\begin{equation}
\label{e:tt.2}
\zeta_t + L_2\zeta + u\cdot\nbl\tht + v\cdot\nbl\zeta
 +w_t\tht_z + w \zeta_z =0.
\end{equation}
 Taking the inner product of \eqref{e:vt} with $u$ and using the boundary
 conditions \eqref{e:bc.v} and \eqref{e:bc.t}, we obtain
\begin{equation}
\label{e:vt.1}
\begin{split}
\hf\dert\|u\|_2^2 +& \|u\|_{V_1}^2\\
=& -\lra{(u\cdot\nbl)v}{u} -\lra{(v\cdot\nbl)u}{u}
  - \lra{w_tv_z}{u} - \lra{w u_z}{u} \\
 & - \lra{\nbl (p_s)_t}{u} - \lra{\int_z^0\nbl\zeta}{u} - \lra{fu^\bot}{u}\\
=& -\lra{(u\cdot\nbl)v}{u} - \lra{w_tv_z}{u} + \lra{\int_z^0\nbl\zeta}{u}
\end{split}
\end{equation}
 where we have used the fact that
\[ \lra{\nbl (p_s)_t}{u} = \lra{fu^\bot}{u}
   = \lra{(v\cdot\nbl)u}{u} + \lra{w u_z}{u} =0. \]
 Taking the inner product of \eqref{e:tt} with $\zeta$ and using \eqref{e:bc.t},
 we obtain
\begin{equation}
\label{e:tt.1}
\begin{split}
 \hf\dert\|\zeta\|_2^2 +& \|\zeta\|_{V_2}^2
 +\mu_2\alf\|\zeta(z=0)\|_2^2\\
=& -\lra{u\cdot\nbl\tht}{\zeta} -\lra{v\cdot\nbl\zeta}{\zeta}
 - \lra{w_t\tht_z}{\zeta} -\lra{w \zeta_z}{\zeta}\\
=& -\lra{u\cdot\nbl\tht}{\zeta} - \lra{w_t\tht_z}{\zeta}
\end{split}
\end{equation}
 where we have used the fact that
\[ \lra{v\cdot\nbl\zeta}{\zeta} + \lra{w \zeta_z}{\zeta} =0. \]
 Following the {\em a priori} estimates of \cite{j:16} using \eqref{e:vt.1} and
 \eqref{e:tt.1}, we obtain
\begin{equation}
\label{e:u.eta}
 \dt{\lf(\|u\|^2+\|\zeta\|^2\rt)} + \|u\|_{V_1}^2 + \|\zeta\|_{V_2}^2
 \lsp g(t) \lf(\|u\|^2+\|\zeta\|^2\rt),
\end{equation}
 where
\[ g(t) :=  1 + \|\nbl v\|_2^4 + \|\nbl\tht\|_2^4 +\|v_z\|_2^2\|\nbl v_z\|_2^2
 + \|\tht_z\|_2^2\|\nbl \tht_z\|_2^2. \]
 Thus, from \eqref{e:u.eta}, we obtain
\[  (v_t, \tht_t)\in L^2( 0,T; V),\ \forall\ T>0. \]
 This finishes proof of Theorem~\ref{t:t}.
\end{prf}

\section{Uniform Boundedness for Higher Regularity}
\label{s:high}

\begin{thrm}
\label{t:hm}
 Suppose the integer $m\gs 3$, $Q\in H^{m-1}(\Omg)$ and
 $(v_0,\tht_0)\in V\cap (H^m(\Omg))^3$. Let $(v,\tht)$ be the unique strong
 sloution of problem \eqref{e:v.n}-\eqref{e:ic.n}. Then,
\[ (v, \tht)\in L^\infty(0, \infty; (H^m(\Omg))^3),\
  (v, \tht) \in L^2(0, T; (H^{m+1}(\Omg)^3),\ \forall\ T\in(0, \infty).\]
 Moreover, there exists a bounded absorbing set of
 $\|(v,\tht)\|_{(H^m(\Omg))^3}$ in $R_+$.
\end{thrm}

\begin{prf}
 The idea of the proof is to use induction and the method of the proof of
 Theorem~\ref{t:h2}. We give an outline of the proof here. To simplify the
 presentation without loss of generality, let us consider the following
 simplified equation:
\begin{equation}
\label{e:v.ps}
 v_t+A_1v= -\Pf\lf[ (v\cdot\nbl)v+wv_z\rt].
 \end{equation}
 Applying $A_1^m$ to \eqref{e:v.ps} and taking inner product with $v$, we get
\begin{equation*}
\begin{split}
 \hf\dt{\|\Lmd_1^mv\|^2} + \|\Lmd_1^{m+1}v\|^2
 =& -\lra{\Lmd_1^{m-1} \Pf\lf[(v\cdot\nbl)v+wv_z\rt]}{\Lmd_1^{m+1}v}\\
 \ls& \hf \|\Lmd_1^{m-1} \Pf\lf[(v\cdot\nbl)v+wv_z\rt]\|^2
  + \hf \|\Lmd_1^{m+1}v\|^2,
\end{split}
\end{equation*}
 that is
\begin{equation}
\label{e:v.ps.1}
\dt{\|\Lmd_1^mv\|^2} + \|\Lmd_1^{m+1}v\|^2 
\ls\|\Lmd_1^{m-1} \Pf\lf[(v\cdot\nbl)v+wv_z\rt]\|^2
\end{equation}
 Now, we can use induction to prove the expected uniform boundedness of
 $\|v\|_{(H^m(\Omg))^2}$ for $m\gs3$. We already have the result valid for
 $m=0,1,2$. Assume it is already valid for $\|v\|_{(H^k(\Omg))^2}$, with
 $0\ls k\ls m-1$. Next, we show uniform boundedness of $\|v\|_{(H^m(\Omg))^2}$.
 Using Theorem~\ref{t:P.d} as in the proof of Theorem~\ref{t:h2}, we have
\[ \|\Lmd_1^{m-1} \Pf\lf[(v\cdot\nbl)v+wv_z\rt]\|^2
 \lsp \|\lf[(v\cdot\nbl)v+wv_z\rt]\|_{(H^{m-1}(\Omg))^2}^2. \]
 To estimate the right-hand side of the above inequality, we just need to
 estimate $L^2$ norm of 
\begin{equation}
\label{e:p3.m-1} \prt_3^{m-1}\lf[(v\cdot\nbl)v+wv_z\rt],
\end{equation}
 since if these terms which contain the highest order derivatives can be well
 controlled, then those terms involving lower order derivatives can be treated
 exactly in the same way. To simplify our presentation, we just need to estimate
\[ \|\Pz^{m-1}[(v\cdot\nbl)v]\|,\ \text{ and }\ \|\nbl^{m-1} (wv_z) \|, \]   
 since all other terms appearing in \eqref{e:p3.m-1} can be treated similarly
 and behave no worse. By Leibniz rule, we have
\[ \|\Pz^{m-1}\lf[(v\cdot\nbl)v\rt]\|
 \ls \sum_{k=0}^{m-1}\|\Pz^kv \cdot\nbl\Pz^{m-1-k}v\|. \]
 By Lemma~\ref{l:ju.2},
\begin{align*}
 \| \Pz^kv \cdot\nbl\Pz^{m-1-k}v\|
\lsp \|v\|_{H^k}^\qt\|v\|_{H^{k+1}}^{\frac{3}{4}}
 \|v\|_{H^{m-k}}^\qt\|v\|_{H^{m-k+1}}^{\frac{3}{4}}.
\end{align*}
 Hence, by induction hypothesis for uniform boundedness of $\|v\|_{(H^{m-1})^2}$,
 we have
\begin{align*}
\|\Pz^{m-1}\lf[(v\cdot\nbl)v\rt]\|
\lsp& \sum_{k=0}^{m-1}\|v\|_{H^k}^\qt\|v\|_{H^{k+1}}^{\frac{3}{4}}
 \|v\|_{H^{m-k}}^\qt\|v\|_{H^{m-k+1}}^{\frac{3}{4}}\\
\lsp& \lf(\max_{0\ls k\ls m-1} \|v\|_{H^k}^\qt\rt)
\sum_{k=0}^{m-1} \|v\|_{H^{k+1}}^{\frac{3}{4}}
 \|v\|_{H^{m-k}}^\qt\|v\|_{H^{m-k+1}}^{\frac{3}{4}}\\
\lsp& \sum_{k=0}^{m-1} \|v\|_{H^{k+1}}^{\frac{3}{4}}
 \|v\|_{H^{m-k}}^\qt\|v\|_{H^{m-k+1}}^{\frac{3}{4}}\\
\lsp& \|v\|_{H^m}^\qt\|v\|_{H^{m+1}}^\qtt + \|v\|_{H^m}^\qtt.
\end{align*}
 Therefore,
\begin{equation}
\label{e:vdv.m-1}
\begin{split}
\|\Pz^{m-1}\lf[(v\cdot\nbl)v\rt]\|^2
 \lsp& \|v\|_{H^m}^\hf\|v\|_{H^{m+1}}^\hft + \|v\|_{H^m}^\hft\\
\lsp& 1 + C_\veps \|v\|_{H^m}^2 + \veps \|v\|_{H^{m+1}}^2.
\end{split}
\end{equation}
 Similarly, we have
\[ \|\nbl^{m-1} (wv_z)\|
 \lsp \sum_{k=0}^{m-1}\lf\| (\nbl^k w)\ox (\nbl^{m-1-k}v_z)\rt\|, \]
 where we use $\ox$ to include all the products occurring here. By \eqref{e:w}
 and Lemma~\ref{l:ju.1}, we get
\begin{equation}
\begin{split}
\label{e:wvz.m-1}
\|\nbl^{m-1} (wv_z)\|^2
 \lsp& \sum_{k=0}^{m-1}
 \|v\|_{H^{k+1}} \|v\|_{H^{k+2}}\|v\|_{H^{m-k}}\|v\|_{H^{m+1-k}}\\
\lsp& \|v\|_{H^m}\|v\|_{H^{m+1}}
\ls C_\veps\|v\|_{H^m}^2 +\veps\|v\|_{H^{m+1}}^2.
\end{split}
\end{equation}
 Therefore, from \eqref{e:v.ps.1} with consideration of\eqref{e:vdv.m-1} and
 \eqref{e:wvz.m-1}, we have for $m\gs 3$,
\begin{equation}
\label{e:v.ps.2}
\dt{\|\Lmd_1^mv\|^2} + \|\Lmd_1^{m+1}v\|^2 
\lsp 1+\|v\|_{H^m}^2.
\end{equation}
 Applying Glonwall lemma first and then uniform Gronwall lemma to
 \eqref{e:v.ps.2} proves uniform boundedness of $\|v\|_{H^m}$ for \eqref{e:v.ps}
 and existence of a bounded absorbing set for $\|v\|_{H^m}$ in $\R_+$. The proof
 of the theorem for problem\eqref{e:v.n}-\eqref{e:ic.n} is similar and is thus
 omitted.

\end{prf}

\end{document}